\documentclass[10pt,reqno]{amsart} 

\usepackage{amssymb,latexsym}
\usepackage{cite} 

\usepackage[height=190mm,width=130mm]{geometry} 

\theoremstyle{plain}
\newtheorem{theorem}{Theorem}
\newtheorem{lemma}{Lemma}
\newtheorem{corollary}{Corollary}

\theoremstyle{definition}

\theoremstyle{remark}

\newcommand{\Z}{\mathbb{Z}}
\newcommand{\C}{\mathbb{C}}
\newcommand{\R}{\mathbb{R}}

\numberwithin{equation}{section} 

\newcommand{\N}{\ensuremath{\mathbb {N}}}

\newcommand{\g}{\ensuremath{\Gamma}}

\newcommand{\ps}{{\raise 1pt\hbox{\tiny (}}}

\newcommand{\pss}{{\raise 1pt\hbox{\tiny [}}}
\newcommand{\pdd}{{\raise 1pt\hbox{\tiny ]}}}
\newcommand{\pd}{{\raise 1pt\hbox{\tiny )}}}

\newcommand{\bs}{{\raise 1pt\hbox{\tiny [}}}
\newcommand{\bd}{{\raise 1pt\hbox{\tiny ]}}}

\def\cross{\mathinner{\mathrel{\raise0.8pt\hbox{$\scriptstyle>$}}
                 \joinrel\mathrel\triangleleft}}

\usepackage{citehack}
\usepackage{amsmath,amsthm}
\usepackage{amsfonts}
\usepackage{amssymb}
\usepackage{longtable}
\usepackage[matrix,arrow,curve]{xy}

\usepackage{lipsum}

\usepackage{amsfonts}
\usepackage{amsmath}
\usepackage{amscd}
\usepackage{amssymb}
\usepackage{amsthm}
\usepackage{bbm}
\usepackage{bm}
\usepackage{enumerate}
\usepackage{enumitem}
\usepackage{hyperref}
\usepackage{latexsym}
\usepackage{mathabx}
\usepackage{mathdots}
\usepackage{mathtools}

\usepackage{amsmath}
\usepackage{amssymb}
\usepackage{amsthm}
\usepackage{bbm}
\usepackage{enumerate}
\usepackage{latexsym}
\usepackage{mathdots}
\usepackage{geometry}

\usepackage{hyperref}  
 
\usepackage{amsmath}
\usepackage{amssymb}
\usepackage{amsthm}
\usepackage{bbm}
\usepackage{enumerate}
\usepackage{latexsym}
\usepackage{mathdots}
\usepackage{geometry}

\def\K{\mathcal{K}}

\usepackage{stackrel}

\newcommand{\be}{\begin{equation}}
\newcommand{\ee}{\end{equation}}

\newcommand{\nn}{\nonumber \\}

\newcommand{\wt}{\mbox{\rm wt}\ }

\newcommand{\nc}{\newcommand}
\nc{\cali}{\mathcal}
\nc{\on}{\operatorname}
\nc{\Wick}{{\mb :}}

\nc{\ddz}{\frac{\partial}{\partial z}}
\nc{\ch}{\mbox{ch}}
\nc{\Oo}{{\cali O}}
\nc{\cond}{|\,}
\nc{\bib}{\bibitem}
\nc{\pone}{\Pro^1}
\nc{\pa}{\partial}
\nc{\arr}{\rightarrow}
\nc{\larr}{\longrightarrow}
\nc{\ket}{\rangle}
\nc{\bra}{\langle}
\nc{\gam}{\bar{\gamma}}
\nc{\q}{\widetilde{Q}}
\nc{\ep}{\epsilon}
\nc{\su}{\widehat{{\mf s}{\mf l}}_2}
\nc{\sw}{{\mf s}{\mf l}}
\nc{\h}{{\mf h}}
\nc{\n}{{\mf n}}
\nc{\ab}{\mf{a}}
\nc{\is}{{\mb i}}
\nc{\js}{{\mb j}}
\nc{\bi}{\bibitem}
\nc{\He}{{\cali H}}
\nc{\inv}{^{-1}}
\nc{\ol}{\overline}
\nc{\wh}{\widehat}
\nc{\dst}{\displaystyle}

\nc{\delt}{\partial_t}
\nc{\ddt}{\frac{\partial}{\partial t}}
\nc{\delx}{\partial_x}
\nc{\mb}{\mathbf}
\nc{\mf}{\mathfrak}

\nc{\mbb}{\mathbb}
\nc{\Ctt}{\C((t))}
\nc{\Ct}{\C[t,t\inv]}

\nc{\ghat}{\wh{\g}}

\nc{\un}{\underline}
\nc{\mc}{\mathcal}
\nc{\BB}{{\mc B}}
\nc{\bb}{{\mf b}}
\nc{\kk}{{\mf k}}
\nc{\frob}{\times}
\nc{\sm}{\setminus}
\nc{\Pp}{{\mathbb P}^1}
\nc{\Aa}{{\mc A}}

\nc{\AutO}{\on{Aut}\Oo}
\nc{\AUTO}{\un{\on{Aut}}\Oo}
\nc{\AUTK}{\un{\on{Aut}}\K}
\nc{\Heout}{\He_{\out}}
\nc{\Hetil}{{\widetilde\He}}
\nc{\wb}{\overline}

\nc{\Res}{\on{Res}}
\nc{\pitil}{\Pi}
\nc{\Ctil}{\wt{C}}
\nc{\auto}{\on{Aut} \Oo}
\nc{\phitil}{\wt{\phi}}
\nc{\gz}{\g_{\vec z}}
\nc{\tensorM}{\bigotimes_{i=1}^N{\mathbb M}_i}
\nc{\tensorW}{\bigotimes_{i=1}^N W_{\nu_i,k}}
\nc{\out}{\on{out}}

\nc{\m}{{\mathfrak m}}

\nc{\gx}{\g^0_{\vec x}}

\nc{\hx}{\He^0_{\vec x}}
\nc{\tensorpi}{\pi_{\nu_1,\ldots,\nu_N}^\kappa}
\nc{\Phizw}{\Phi_{\vec w}({\vec z})}
\nc{\Pro}{{\mathbb P}}

\nc{\De}{\overline{d}}

\nc{\us}{\underset}

\nc{\Ll}{\mc L}
\nc{\dR}{\on{dR}}

\nc{\T}{{\mc T}}

\nc{\Xn}{\overset{\circ}X{}^n} \nc{\Dn}{\overset{\circ}D{}^n}
\nc{\Dxn}{\overset{\circ}D{}^n_x} \nc{\varphitil}{\wt{\varphi}}

\nc{\lf}{{\mf l}}
\nc{\GL}{{}^L G}
\nc{\Vir}{\on{Vir}}

\begin{document}
\title[The hierarchies of identities and closed products]  
{The hierarchies of identities and closed products for multiple complexes}   

\author{Daniel Levin${}^\&$}
\address{ ${}^\&$ Mathematical Institute \\ University of Oxford \\ Andrew Wiles Building \\
 Radcliffe Observatory Quarter (550) \\ Woodstock Road \\ Oxford \\ OX2 6GG \\ United Kingdom 
(Currently working in Israel)}

\author{Alexander Zuevsky${}^*$}
\address{ ${}^*$ Institute of Mathematics \\ Czech Academy of Sciences\\ Zitna 25, 11567 \\ Prague\\ Czech Republic }

\email{levindanie@gmail.com}
\email{zuevsky@yahoo.com}

\begin{abstract}
We consider infinite $\Z_\Z$-index complexes  
$\mathcal C$  
of spaces with elements depending on a number of parameters, 
complete with respect to 
a linear associative regular inseparable multilinear product.  
The existence of nets of vanishing ideals   
 of orders of and powers of differentials 
is assumed for subspaces of $\mathcal C$-spaces. 
In the polynomial case of orders and powers of the differentials, 
we derive the hierarchies of differential identities 
and closed multiple products. 
We prove that a set of maximal orders and powers for differentials, 
 differential conditions,       
together with coherence conditions 
 on indices of a complex $\mathcal C$ elements 
generate families of multi-graded differential algebras. 
AMS Classification: 53C12, 57R20, 17B69 
\end{abstract}

\keywords{Differential identities for multiple complexes; 
combinatorics of closed products; 
 graded differential algebras} 
\vskip12pt  

\maketitle

The authors state that: 

1.) The paper does not contain any potential conflicts of interests. 

2.) The paper does not use any datasets. No dataset were generated during and/or analysed 
during the current study. 

3.) The paper includes all data generated or analysed during this study. 

4.) Data sharing is not applicable to this article as no datasets were generated or analysed during the current study.

5.) The data of the paper can be shared openly.  

6.) No AI was used to write this paper. 
\section{Introduction} 
\label{valued}
The main purpose of this paper is derive 
hierarchies of differential identities  
and closed products following 
from sets of natural conditions on orders and powers of 
differentials applied to elements of a multiple complex.  
In contrast to ordinary wedge-product case for differential forms 
 reflected in \cite{Ghys, Ko} 
we work with the universal enveloping algebra 
constituted by $\N$-valued powers of elements of complexes,  
in particular, given by differentials applied to such elements. 
Note that in general we do not specify commutation relations 
for elements of $\mathcal C$.  
Nevertheless, the structure given in the paper
 reproduces the structure of a differential algebra. 
The differential identities for elements of 
complexes with multiple indices
 endowed with regular associative products  
is an important way to study various algebraic and geometric structures. 
In particular, they are extremely useful to find closed products in 
 cohomology classes computations  
of invariants associated to a multiple complex.
 Differential conditions applied to elements of a 
multiple complex $\mathcal C$ provide 
a system of 
multiplication rules and form the resulting algebra. 
The full algebra associated to all possible choices of differential relations
is quite huge. 
For differential forms considered on smooth manifolds,  
the Frobenius theorem for a distribution 
 leads to orthogonality conditions on forms.      
Motivated by the notion of integrability 
for differential forms on foliated manifolds 
\cite{Ghys, Ko},    
 we study systems of differential identities 
with respect to a multiple product of elements 
of a multiple-index complex. 
 We then show that differential conditions  
endow a multiple complex with the structure of 
a graded differential algebra. 
Assuming existence of a net of differential power/order vanishing ideals 
 in subspaces of a multiple complex $\mathcal C$,   
and requiring natural orthogonality conditions 
for completions of a set of $\mathcal C$ elements,  
 the hierarchies of differential conditions arise in terms of closed products. 

Ideologically, closed products 
containing powers of multiple action of mixed differentials 
represent a geometrical idea of ``codimension`` one products. 
Namely, there are two possibility. 
First, a product contains a combination of 
differentials which vanishes due to the critical orders or powers of 
 differentials  
as the result of application of differentials. 
Second, an initial candidate to closed products, 
does not contain maximal orders or powers 
but after an application of a differential identity, 
the total maximal powers of elements appear, and  
a product can be reduced to the first case.  
The hierarchies of differential identities are of non-trivial nature 
since we assume that  maximal orders of differentials 
as well as maximal powers of $\mathcal C$-elements  
depend on elements themselves 
(see explanations in Section \ref{setup}). 
Certain non-trivial examples of multiple product identities 
are given in Sections \ref{porgo}--\ref{mima}.   
We provide examples of closed products 
with respect to multiple products of 
differentials $d_a$ and $d_{\overline{a}}=\overline{d}_a$ 

The hierarchies we derive from the conditions on a complex
are useful in the theory of continual Lie algebras 
\cite{sv1, sv2, sv3, vershik, ass}
and in the theory of completely integrable \cite{arnold}
and exactly solvable \cite{ls, razsav, bakas, FMc}
 dynamical systems. 
In particular, similar to invariants associated to foliations, 
one is able to use the identities from a hierarchy 
to prove integrability as well as to find invariants 
of corresponding dynamical systems. 
It is important to mention that the hierarchies of identities
and closed products  
(i.e., products annihilating by 
a single differential) 
constitute the tools for direct computation and classification of 
cohomology invariants of the corresponding complex. 
Therefore, we are interested in generating all possible 
closed products. 
We will present such a classification in a forthcoming paper. 
The next step in finding invariants is to prove their 
independence with respect to replacements of a complex elements. 
The classification problem 
of cohomology invariants associated to a complex 
endowed with a multiple product 
will be treated in a separate paper. 

As for possible applications of the material presented in this paper,
 we would like to mention 
computations of higher cohomology for grading-restricted vertex algebras \cite{Huang},  
search for more complicated 
cohomology invariants, and applications in differential geometry and algebraic topology.  
It would be interesting to study possible applications of invariants we constructed
 to cohomology of manifolds. 
In differential geometry there exist various 
approaches to the construction of cohomology classes 
(cf., in particular, \cite{Losik}).
 We hope to use these techniques to derive counterparts in 
the cohomology theory of vertex algebras.
The results proven in this paper are also 
useful in computations of cohomology of foliations \cite{BG, BGG}.  
Finally, results of this paper will be useful for computations 
\cite{zub2016, Zubkov-2508, Zubkov-Wu}  
of topological 
invariants in quantum field theory via cyclic cohomology methods \cite{Connes}.    
In addition to that, the constructions of this paper 
will be useful in descriptions of 
 relativistic quantum field theory and 
fermionic superfluids \cite{z7}, 
the theory of topological invariants in
high energy physics \cite{z4, z5, z6}, 
chiral separation effect \cite{z3},  
various topological defects dominates dynamics \cite{z8},  
 non-renormalization by interactions of Integer quantum Hall effect \cite{z9}, and 
Wigner-Weyl calculus \cite{z1, z2}. 
\section{Multiple complexes and the main result}
\label{setup}
Introduce a system of families of multiple horizontal 
and vertical complexes   
$\mathcal C=
\left(C\left(\Theta^{\bm n}_{\bm m}\right), 
 d^{\bm n'}_{\bm m'}, {\overline{d}}^{\bm n''}_{\bm m''} \right)$
 with ${\bm n}\in \Z_\Z$ and ${\bm m}\in \Z_\Z$, 
i.e., an infinite number of up and down $\Z$-valued indices, 
and ${\bm n}={\bm n}' \cup {\bm n''}$, ${\bm m}={\bm m}' \cup {\bm m''}$.    
These indices correspond to increasing and decreasing indices of parameters  
$\Theta^{\bm n}_{\bm m}$ 
for spaces $C\left(\Theta^{\bm n}_{\bm m}\right)$ 
under actions of $\mathcal C$-differentials  
correspondingly. 
Here $d^{\bm n'}_{\bm m'}$ denotes a family of differentials 
for each pair of entries $n_i'$ of ${\bm n'}$ and $m_i'$ of ${\bm m'}$, $i \in \Z$.  
The action of
a differential 
$d^{n_i'}_{m_i'}: C\left(\Theta^{\bm n}_{\bm m}\right) 
\to C\left(\Theta^{\ldots, n_i' + 1, \ldots,  }_{\ldots, m_i' - 1, \ldots}\right)$,
 $i \in \Z$, 
for ${\ldots, n_i'+1, \ldots}$ and ${\ldots, m_i'-1, \ldots}$
 means that all entries  
of the multiindices ${\bm n'}$ and ${\bm {m'-1}}$ 
 remain the same except for   
$n_i' \mapsto n_i'+1$ and $m'_i \mapsto m_i-1$,   
that increases or decreases by one 
correspondingly. 
Other indices $\bm n''$ and $\bm m''$ remain unchanged.
 
The family of vertical differentials  
$\overline{d}^{\bm n''}_{\bm m''}: 
C\left(\Theta^{\bm n}_{\bm m}\right) 
\to C\left(\Theta^{\ldots, n_j'' + 1, \ldots}_{\ldots, m_j'' - 1, \ldots}\right)$, 
$j \in \Z$, 
change $n''_j$ and $m''_j$ indices similarly, 
while the indices $\bm n'$ and $\bm m'$ remain unchanged.
A subset of horizontal and vertical complexes 
corresponding to subsets of upper indices in ${\bm n}$ 
and lower indices in ${\bm m}$ 
are supposed to be chain, cochain, or chain-cochain 
corresponding to 
the differentials $\left(d\right)^{\bm n}_{\bm m}$  
 and $\overline{d}^{\bm n}_{\bm m}$. 
When we write $\chi \in \mathcal C$ that means that an element $\chi$ belongs to 
a subspace of $\mathcal C$.  
The $i\in \Z$, $j\in \Z$ and $k \in \Z$-th slice of the full diagram of  
 a complex $\mathcal C$ is described by the following diagram 
\begin{eqnarray*}
    \quad   \vdots 
\qquad \qquad \qquad \qquad \qquad \qquad   \vdots
\qquad \qquad \qquad \qquad \qquad \qquad \qquad \quad 
 &&
\nn
    \downarrow 
\overline{d}^{\ldots, n_i'-1, \ldots, n_j''-1, \ldots }_{\ldots, m'_i,  \ldots, m''_j, \ldots}
\qquad  \; \qquad  \downarrow  
 \overline{d}^{\ldots, n_i-1, \ldots, n_j-1, \ldots}_{\ldots, m_i'-1, \ldots, m''_j, \ldots }
 \qquad \qquad \qquad \quad &&
\nn
 \ldots \stackrel{d^{\ldots, n_i-1, \ldots, n_j, \ldots}_{\ldots, m_k+1, \ldots} }
 {\longrightarrow} 
C^{\ldots, n_i, \ldots, n_j, \ldots}_{\ldots, m_k, \ldots}  
   \stackrel{d^{\ldots, n_i, \ldots, n_j, \ldots}_{\ldots, m_k, \ldots} }
 {\longrightarrow }
   C^{\ldots, n_i+1, \ldots, n_j, \ldots}_{\ldots, m_k-1, \ldots} 
  \stackrel{d^{\ldots, n_i+1, \ldots, n_j, \ldots}_{\ldots, m_k-1, \ldots} } 
{\longrightarrow}
 \cdots \qquad \quad && 
\nn
 \qquad \qquad \quad  \downarrow   
\overline{d}^{\ldots, n_i, \ldots, n_j, \ldots}_{\ldots, m_k, \ldots} 
\qquad  \quad \quad \qquad 
\downarrow  \overline{d}^{\ldots, n_i+1, \ldots. n_j, \ldots }_{\ldots, m_k-1, \ldots}
 \qquad \qquad  \qquad \qquad && 
\nn
 \ldots \stackrel{d^{\ldots, n_i-1, \ldots, n_j+1, \ldots }_{\ldots, m_k+1, \ldots} }
 {\longrightarrow} 
C^{\ldots, n_i, \ldots, n_j+1, \ldots }_{\ldots, m_k, \ldots} 
    \stackrel{d^{\ldots, n_i, \ldots, n_j+1, \ldots}_{\ldots, m_k, \ldots} }
 {\longrightarrow} 
 C^{\ldots, n_i+1, \ldots, n_j+1, \ldots}_{\ldots, m_k-1, \ldots}
  \stackrel{d^{\ldots, n_i+1, \ldots, n_j+1, \ldots}_{\ldots, m_k-1, \ldots}}   
{\longrightarrow} 
 \cdots  &&
\nn
 \qquad \qquad \quad    \downarrow 
 \overline{d}^{\ldots, n_i, \ldots, n_j+1, \ldots}_{\ldots, m_k, \ldots} 
\qquad  \quad \quad \qquad  \downarrow  
\overline{d}^{\ldots, n_i+1, \ldots, n_j+1, \ldots}_{\ldots, m_k-1, \ldots}
 \qquad \qquad  \qquad \qquad  &&
\nn
\quad   \vdots 
\qquad \qquad \qquad \qquad \qquad \qquad   \vdots
\qquad \qquad \qquad \qquad \qquad \qquad \qquad \quad && 
\end{eqnarray*} 
Introduce the following notations.
Since indices of differentials and indices of $\mathcal C$-subspaces 
are taken to be always coherent, let us denote  
$d_a =\left(d_a\right)^{\bm n}_{\bm m}$, 
where the index $a$ denotes either of the differentials $d$ or $\overline{d}$, 
$\overline{a}$ denotes the opposite choice, 
and $d_{\overline{a}}=\overline{d}_{a}$. 
A combination of differentials is notated by $d^{r_1}_{a_1} \ldots d^{r_k}_{a_k}$,   
where $r_1$, $\ldots$, $r_k$ are orders of corresponding differentials.  
We call $q$ of $(\phi)^q$ the power of an element $\phi$. 
For all possible combinations of $l$ elements of $\mathcal C$  
we assume the existence of a formal multi-linear 
 associative inseparable $\cdot_l$-product 
\begin{equation}
\label{tugo} 
 \cdot_l = \cdot_{i=1}^l: \times_{i=1}^l C^{ {\bm n}_i }_{ {\bm m}_i }   
\to C^{ \widetilde{\bm n} }_{ \widetilde{\bm m} },     
\end{equation}  
with $\widetilde{\bm n}(\phi)=\sum\limits_{i=1}^l {\bm n}_i(\phi_i)$,
$\widetilde{\bm m}(\phi)=\sum\limits_{i=1}^l {\bm m}_i(\phi_i)$
where the summations are component-vise.   
We set that the result of a $\cdot_l$-product 
which is an element of a $\mathcal C$-space, 
is regular in the domains of definition
 of all parameters $\Theta^{\bm n}_{\bm m}$ defining the product. 
We assume that a $\cdot_l$-product is inseparable in 
pairs of two elements of $\mathcal C$ in general. 
Let us introduce the notations. 
In the arguments of a multiple product we denote by 
$(\phi)_i$ a $\mathcal C$-element placed at the $i$-th position.  
The notation $\widehat{s}$ means that the value $s$ of the index 
is omitted.   
We assume the completeness of a complex $\mathcal C$-spaces  
with respect to the multiple product \eqref{tugo}. 
That means that for $\phi$, $\psi \in \mathcal C$, 
placed in the product at the $i$-th and $j$-th positions, 
  from $0=\cdot_l(\ldots, (\phi)_i, \ldots, (\psi)_j, \ldots)$
it follows that there exists a $\mathcal C$-element $\chi$, 
such that 
$\psi=\cdot_l(\ldots, (\phi)_{i'}, \ldots, (\chi)_k, \ldots)$, 
with some $i'$ and $k$. 
For a $\cdot_l$-product we set $\cdot_l(\ldots, 0, \ldots)=0$. 
One can also set that 
$\cdot_l(\ldots, {\rm Id}_\mathcal C, \ldots)=\cdot_{l-1}(\ldots, 0, \ldots)=0$,  
where ${\rm Id}_\mathcal C$ 
is the identity element in the corresponding $\mathcal C$-subspace.   
In certain cases, a product \eqref{tugo} may have coincidences of   
parameters of multiplied elements of a multiple complex $\mathcal C$ spaces.  
The result of a product may not allow such coincidences leading to 
possible overcounting of the number of parameters. 
In order avoid such a possibility, 
we take into account one coinciding elements/parameters only
 in the result of a product.  
For certain $\cdot_l$-products of $\mathcal C$-elements depending on 
parameters, 
(e.g., for a product of vertex operator algebra complexes \cite{Huang}),  
 it is needed to exclude 
a number of coinciding arguments 
and count them only once in the resulting product. 
Examples of such products can be found in \cite{Y, Gui} and others. 
For an element $\phi_i \in C^{\bm n}_{\bm m}$ 
in a particular $\cdot_l$-product  
let $\bm r$ and $\bm t$ 
be the numbers of common parameters 
corresponding to upper and lower indices for $\phi$ 
with other elements 
in a product.   
In that case, the conditions on indices for a resulting element 
of a $\cdot_l$-product $\phi$ are 
${\bm n}(\phi) = \sum\limits_{i=1}^l {\bm n}_i(\phi_i)-{\bm r}_i(\phi_i)$, 
and ${\bm m}(\phi)= \sum\limits_{i=1}^l {\bm n}_i(\phi_i)-{\bm t}_i(\phi_i)$.  
The associativity of a $\cdot_l$-product means  
$(a\cdot_l \ldots \cdot_{l'} b) \cdot_{l''} \ldots \cdot_{l'''}c
=a \cdot_l \ldots \cdot_{l'} (b \cdot_{l''} \ldots \cdot_{l'''} c)$. 
 for all elements $a$, $b$, $c \in \mathcal C$.  
Note that $l$ can be taken infinite if we assure that a product 
$\times_{i \ge 1}$ is converging with infinite $l$.  
We call the product formal since elements of a complex $\mathcal C$
spaces can be formal (in particular, geometric) objects. 
Then the result of a $\cdot_l$-product is 
 a superposition 
of formal objects (e.g., Riemann surfaces \cite{Y}).    
In that case, convergence of a product means that  
 the corresponding superposition leads to a well-defined formal object. 

In this paper, 
all the constructions are independent of actual commutation relations  
 for elements for elements inside a $\cdot_l$-product. 
In addition to all that above we assume that  
in a multiple complex $\mathcal C$ spaces   
there exist a net of subspaces constituting 
  vanishing products of exactly $q$ its elements.   
 That we call an order $q$ ideal $\mathcal I(q) \subset \mathcal C$.  
We define a distributed ideal 
$\mathcal I(q)$ of order $q \le l$ as   
a union of subsets of $\mathcal C$ such that 
for any set of $q$ elements 
$\theta_1$, $\ldots$, $\theta_q \in \mathcal I(q)$, 
distributed in a product, 
$\cdot_l (\ldots$, $\theta_1$, $\ldots$,  $\theta_2$,  
$\ldots$, $\theta_{q-1}$, \ldots, $\theta_q$, $\ldots )=0$.   
A product \eqref{tugo} vanishes if at least one entry of its arguments 
belonging to $\mathcal C$ is zero. 
Now, let us explain how we understand powers of 
 an element $\phi \in \mathcal C$. 
We denote 
$(\ldots, (\phi)^r, \ldots )
=(\ldots, (\phi)_{j_1}, \ldots, (\phi)_{j_2}, \ldots, (\phi)_{j_r}, \ldots)$  
for $\phi$ placed at the positions $(j_1, \ldots, j_r)$ with  
 $r \le q(\phi)$ where $q(\phi)$ is the maximal power of an element $\phi$, 
i.e., for $r \ge q(\phi)$,  $(\phi)^r=0$. 

 It is taken that a product \eqref{tugo} 
is coherent in indices with respect to all differentials 
of a complex $\mathcal C$.  
That means that for every set of 
$\phi_i \in C\left(\Theta^{\bm {n_i}}_{\bm {m_i}}\right)$, $1 \le i \le l$,  
 ${\bm n}(\phi)=\sum\limits_{i=1}^l {\bm {n_i}}(\phi_i)$,  
 and ${\bm m}(\phi)=\sum\limits_{i=1}^l {\bm {m_i}}(\phi_i)$,  
and the resulting 
$\phi \in C\left(\Theta^{\bm n}_{\bm m}\right)$, 
$\phi=\cdot_l(\phi_1, \ldots, \phi_l)$,   
 the action of corresponding differential $d$   
satisfies Leibniz formula 
\begin{equation}
\label{leibniz}
d_a\phi = d_a \cdot_l    
\left(\ldots, (\phi)_i, \ldots \right)=  
\sum\limits_{i=1}^l \cdot_l 
\left(\ldots, \left(d_a \phi\right)_i, \ldots \right). 
\end{equation}
Note that in general, being applied 
to the result of a $\cdot_l$-product,  
 Leibniz formula  
represents a map $\mathcal C\to \times \mathcal C$.  
Therefore a result of Leibniz formula application vanishes  
if in each element of the sum \eqref{leibniz} an element is zero. 

For a set of elements $\phi_k$, $1\le k \le l$, 
in a $\cdot_l$-product 
\begin{equation}
\label{mama}
\cdot_l(\Phi_1, \phi_1, \Phi_2, \ldots, 
\Phi_i, \phi_i, \Phi_{i+1}, \ldots, \Phi_k, \phi_k, \Phi_{k+1}),  
\end{equation}
we call the union of sets of $\mathcal C$-elements
 $\left\{ \Phi_j \right\}$, $1 \le j\le k+1$, 
the completion of a set of elements $\phi_k$ in a $\cdot_l$-product. 
Due to the associativity of a $\cdot_l$-product 
 we can view each set $\Phi_j$ as a single element.
The completion \eqref{mama} for a set of elements $\phi_i$, $1 \le i \le k$,  
with respect to a differential $d_a$   
 is called closed if 
\begin{equation}
\label{papa}
\sum\limits_{j=1}^{p+1} d_{a, j}.\cdot_l(\Phi_1, \phi_1, \Phi_2, \ldots, 
\Phi_i, \phi_i, \Phi_{i+1}, \ldots, \Phi_k, \phi_k, \Phi_{k+1})=0,   
\end{equation}
where $d_{a, j}$, $j \le k+1$, acts on $\Phi_j$-elements of the completion only.  
In what follows we will always assume
that all completions in $\cdot_l$-products
 are closed with respect to corresponding differentials.
Thus, we will skip the completion elements $\Phi_i, \ldots, \Phi_{k+1}$ 
of $\mathcal C$-elements $\phi_i$, $1 \le i \le k$,    
so that a $\cdot_l$-product 
will be denoted as $(\ldots, \phi_1, \ldots, \phi_k \ldots)$.  
Due to the associativity property of the $\cdot_l$-product mention above, 
we will skip the notation $\cdot_l$ and denote all 
completion closed products for various $l$  
as $(\ldots, \ldots, \ldots)$. 
Since all actions of $\mathcal C$ differentials are coherent   
with the indexing of $\mathcal C$-spaces, 
 we will skip also upper and lower  
indices from the notations of $d^{n, \kappa}_m$.  

In general, an arbitrary element $\phi \in \mathcal C$ is  
characterizes by its maximal power 
$q(\phi) \in \Z_+\cup \left\{0\right\}$, i.e.,    
$\left( \ldots, (\phi)_{i_1}, \ldots, (\phi)_{i_{q(\phi)}}, \ldots \right)$ $=0$
 with $q$ identical elements $\phi$ 
distributed in a product \eqref{tugo}.
Next, according to values of maximal orders and powers
 of a collection of differentials $d_{r_i}^{p_i}$, $1 \le i \le k$,   
a $\mathcal C$-element 
$\phi$ satisfies certain conditions with respect to 
subsequent actions of differentials $d_{r_i}^{p_i}$ on $\phi$.   

One can also introduce 
the rules
(possible commutation relations, chain-cochain property, 
maximal orders and powers) 
 for actions of $\mathcal C$-differentials $d_b$, 
with $b=a$ or $\overline{a}$, i.e., 
\begin{eqnarray}
\label{roda}
 d_a d_b= A_{a, b} \; d_b d_a, 
\end{eqnarray}
with some $A_{a, b}$.  
Note that $A_{a, b}$ can be either 
 a complex number (in particular, zero), a map $\mathcal C\to \mathcal C$, 
 or undefined at all. 
We assume that for an element $\phi \in \mathcal C$ 
 might exists the maximal order $q(\phi) \in \R$ depending on $\phi$, 
 such that 
for $r \in \R$,  $r \ge q(\phi)$, $(\phi)^r=0$. 
Similarly, for both differentials $d_b$, $b=a$ or $\overline{a}$, 
 and any element $\phi \in \mathcal C$, 
might exists $p(\phi) \in \R$, depending on $\phi$, 
such that for any $s \ge p(\phi)$,   
$d^s\phi=0$.  

We assume always 
that when we apply a differential $d_a$ to a product $(\gamma)$,  
then there exists $p(\gamma) \in \R$ such that for all $s \ge p(\gamma)$ 
$d_a^s(\gamma)=0$, 
 even if we do not know exact value of $p(\gamma)$. 
In deriving identities for $\mathcal C$-elements,  
 we should always keep in mind 
that the power of a collection of differentials 
applied to a $\mathcal C$-element  
can overcome corresponding maximal value and, 
therefore, a result of such applications 
vanish at some iteration.  
Note again that an element $d_a^{p_i}(\gamma)$, 
$0 \le p_i < p(d_a (\gamma))$ is characterized by 
$q(d_a^{p_i}(\gamma))$ such that 
$(d_a^{p_i}(\gamma))^{q(d_a^{p_i}(\gamma))}=0$. 
For $p_i=0$,  $q(d_a^{p_i}(\gamma))=q(\gamma)$. 
The property of having maximal orders of differentials 
 and powers of elements
 can be seen as a redefinition  
of an initial complex $\mathcal C$ to resulting complexes   
with differentials defined by all possible partitions of all   
 $\mathcal C$ differentials into pairs of differentials  
$d_\sigma$ and $\overline{d}_\sigma$ 
where $\sigma$ is a particular element of a partition. 

Let us introduce some further notations. 
For a set of double indices 
${\bm J}=(J_1, \ldots, J_n)$, $J_i=\left({r_i \atop a_i}\right)$, 
define 
\begin{equation}
\label{pirdos}
\mathcal D_{\bm J} \phi 
=d_{a_1}^{r_1}\ldots d_{a_n}^{r_n} \phi. 
\end{equation} 
 The indices $a_i$ denote the type of differentials. 
The indices $r_j$ are orders of differentials 
satisfying the recurrence conditions
\begin{equation}
\label{pirdo}
 r_j < p\left(d^{r_{j-1}}_{a_{j-1}} \ldots d^{r_1}_{a_1}\phi\right),
\end{equation} 
with the corresponding maximal order for the differential $d_{a_j}^{r_j}$ 
acting on a $\mathcal C$-element 
$d^{r_{j-1}}_{a_{j-1}}$ $\ldots$ $d^{r_1}_{a_1}\phi$, 
$1 \le j \le n_i$. 
Note, that due to the property \eqref{leibniz} and the definition of 
 powers of an element distributed in a $\cdot_l$-product,  
we do not take into account conditions on $\phi \in \mathcal C$ of the form 
$(d_{a_n}^{r_n}(\ldots(d_{a_2}^{r_2}(d_{a_1}^{r_1} \phi)^{q_1}))^{q_2}\ldots)^{q_n}$ 
since, according to \eqref{leibniz} 
 it is equivalent to a finite sum of powers of orders of differentials. 
Note that the definitions of maximal orders and powers above 
characterize elements of a complex $\mathcal C$, and 
 are given in the form of 
differential/orthogonality 
relations \eqref{groby} of Section \ref{bardos}.     
In what follows,  
we deal with multiple products 
 completion of a number of $\mathcal C$-elements 
closed with respect to corresponding differentials. 

Let us now formulate the main result of this paper: 
\begin{theorem}
\label{buldo}
The conditions \eqref{roda} together with a set of maximal orders and powers of 
differentials for a multiple complex $\mathcal C$ 
result in a hierarchy of closed products in terms of differential identities 
on elements of $\mathcal C$ 
 given by the general formula 
\begin{eqnarray}
\label{hier}
0= \sum\limits_{{\bm J}_1, \ldots, {\bm J}_k}
 \left(\ldots, \left(\mathcal D_{{\bm J}_1} \phi_1\right)^{q_1},  
\ldots, \left(\mathcal D_{{\bm J}_k} \phi_k\right)^{q_k}, \ldots \right), 
\end{eqnarray}
with $1 \le q_i < q\left(\mathcal D_{ {\bm J}_i} \right)$,  $1 \le i \le k$,  
where $q\left(\mathcal D_{{\bm J}_i}\right)$ 
 are the maximal powers for the corresponding differentials. 
The differentials $\mathcal D_{\bm J_{s, i}}$ are 
determined by all vanishing products closed with respect to 
single differentials $d_a$. 
\end{theorem}
Now let us give a proof of Theorem \ref{buldo}. 
\begin{proof}
The idea of the proof is related to geometry. 
Namely, we show that the identities \eqref{hier} 
 are defined by vanishing products  
closed with respect a differential $d_a$.  
Geometrically, this corresponds to a separation of 
``codimension one`` differential forms \cite{Ko} expressed
 in terms of multiple products. 
The main idea to generate differential identities for 
multiple products of elements of $\mathcal C$
is to pick products vanishing due to maximal powers  
of differentials action.  
It is easy to see that 
in order to find a differential identity, 
one has to pick elements given by multiple products 
which would give ``almost`` the rule \eqref{pirdo}, 
and would vanish 
 under actions of the differentials 
$d_a$ or $d_{\overline{a}}$.  
Thus, the set of closed products consists 
of such ``almost`` rule \eqref{pirdo}  
elements having ``almost`` 
maximal powers.  
There are two general ways how to obtain an identity
starting from a vanishing product closed under a differential. 
Let us start with a product containing $k$ entries of 
$q_i(\phi_i)$-powers  
of arbitrary $\mathcal C$-elements $\phi_i$, i.e., 
$\left(\ldots, (\phi_i)^{q_i} ,
 \ldots, (\phi_j)^{q_j}, \ldots \right)$,  
 for $0 \le q_i < q(\phi_i)$, with 
the maximal power $q(\phi_i)$ of elements $\phi_i$  
for $1 \le i\ne j \le k$, and,
$q_j = q(\phi_j)$. That product vanishes.  
Recall that for elements $\phi_i$ of $(\phi_i)^{q_i}$     
may be distributed as fingle $\phi_i$ entries among arguments 
of the product. 
Acting by a differential $d_a$
 according to the formula \eqref{leibniz}   
we obtain two types of summands.  
The action of a differential on the distributed power of 
$(\phi_j)^{q(\phi_j)}$ gives $(d_a\phi_j)$ as one entry,  
and remaining distributed entries $(\phi_j)^{q(\phi_j)-1}$.  
As we see, since the maximal power $q(\phi_j)-1$ is
dropped now by one, that particular summand does not vanish. 
There are all together $q(\phi_j)$ non-vanishing summands 
of this type. 
The second type of summands with a differential acting 
on $(\phi_i)^{q_i}$, $i \ne j$,  
 contains the maximal power of $(\phi_j)^{q(\phi_j)}$,  
and, therefore, vanishes. 
Another possibility for a vanishing product to be closed 
is to collect the maximal power of an element 
at some iteration of differentials application.  
We thus obtain the first identity of the hierarchy:
\begin{equation}
\label{pordo}
0=\sum\limits_{s=1}^{q(\phi_j)} \left(\ldots, 
 \phi_i^{q_i} , \ldots, (d_a\phi_j)_s, \ldots, 
 (\phi_j)_{\widehat{s}}^{q(\phi_j)-1}, \ldots \right). 
\end{equation}
This identity relates powers of elements $\phi_i$, $1 \le i \le k$
with the differential $d_a\phi_j$.  
To obtain further identities of that branch of the hierarchy, 
 we apply differentials $d_a$ or $d_{\overline{a}}$ to \eqref{pordo}, 
and take into account \eqref{roda} to derive the next identity in that branch of   
the hierarchy.  
Continuing the process, 
we decrease the powers of elements $(\phi_i)^{q_i}$,   
and increase powers of differentials  
$d^{r_1}_{a_1} \ldots d^{r_n}_{a_n}\phi_i^{q(\phi_i)}$.  
Note that by assumption that the maximal order of 
each differential (i.e., the power when it vanishes)
 applied depends on the element it acts on. 
Therefore, the final form of an identity in the hierarchy 
depends on a sequence of elements 
$d^{r_t}_{a_t} (d^{r_{t-1}}_{a_{t-1}} \ldots (d^{r_1}_{a_1}\phi_{i'})))$, 
$d^{r_{t-1}}_{a_{t-1}}( \ldots (d^{r_1}_{a_1}\phi_{i'}))$,    
$\ldots$, $d^{r_1}_{a_1}\phi_{i'}$, 
where $r_p(\phi_p')$ are lower than the maximal orders 
for corresponding differentials. 
A sequence of identities stops when at least one order 
of differentials reaches its maximal value in each summand. 
The whole hierarchy depends 
on a set of initial elements $\phi_i$. 
That elements may initially contain orders of  
differentials of some powers. 
\end{proof}
To finish this Section we formulate 
\begin{corollary}
The hierarchies \eqref{hier} result in the set of closed multiple products. 
\end{corollary} 
Indeed, each a differential identity $(\chi)=0$ from the system \eqref{hier}
can be ``integrated`` with respect the differentials $d$ or $\overline{d}$ 
to find a multiple product $(\Gamma)$ such that $d_a(\Gamma)=(\chi)=0$.  
When we say integrated we mean that for an element 
$\gamma \in \mathcal C$ we find an element $\gamma'$ 
such that $\gamma'=d_a \gamma$.   
In addition to that relations specified by 
the identities \eqref{hier} can be used 
in order to derive closed products which do not follow
explicitly from the integration of identities.  
Let us underline that we work with a multiple complex spaces as 
grading subspaces of an associative algebra.  
In particular, the multiple product $\cdot_l$ can be introduces for 
geometric objects of general kinds.  
Thus, we keep the constructions 
of this paper independent of commutation relations. 
In the next three Sections we will provide instructive examples
 of identities and closed products 
of the first order of differentials actions.  
The full description and classification 
 of differential identities and closed products 
at higher orders and powers of differentials action  
 will be given in a forthcoming paper. 
The identities and closed products given in Sections \eqref{porgo}--\ref{mima}
 generalize to symmetrized versions with respect to 
permutations of indices, positions of $\mathcal C$-elements, 
and orders and powers of differentials. 
\section{Examples of identities for multiple products} 
\label{porgo}
In this Section we illustrate the proof of Theorem \ref{buldo}  
by complexity-growing examples. 
In this and next two Sections we always assume that all completions 
are closed with respect to corresponding differentials. 
\subsection{Two differentials, and the total maximal power two
 for two kinds of elements} 
 Consider two elements $(\phi, \psi) \in \mathcal I(2)$. 
Recall that the notation $(\phi, \psi) \in \mathcal I(2)$ means 
that a pair of $\phi$ and $\psi$ (but not their differentials) placed 
in a product makes it vanishing.  
 For any two positions   
$1 \le i \ne j \le l$ in a multiple product 
\begin{eqnarray*}
0&=&d_a\left(\ldots, (\phi)_i, \ldots, (\psi)_j, \ldots\right)  
\nn
&=& (\ldots, (d_a\phi)_i, \ldots, (\psi)_j, \ldots)  
+  (\ldots, (\phi)_i,  \ldots, \left(d_a\psi\right)_j, \ldots),     
\end{eqnarray*}
which is not a trivial identity   
for $p(\phi)$, $p(\psi)>1$. 
Recall that the notation $(\phi)_i$ means that an element $\phi$ is placed 
at $i$-th position in a product. 
Therefore, 
\begin{eqnarray}
\label{perdozo}
 (\ldots, (d_a\phi)_i, \ldots, (\psi)_j, \ldots) 
 =-(\ldots, (\phi)_i, \ldots, \left(d_a\psi\right)_j, \ldots), 
\end{eqnarray}
i.e., one can transfer a differential between 
 $\mathcal I(2)$-elements 
while changing the sign. 
In particular, 
 with $d_b d_a=0$, $b=a$ or $b=\overline{a}$, 
and $\overline{p}(\phi)>1$, $\overline{p}(\phi)>1$,  
 we get   
\begin{eqnarray}
\label{pozor}
 &&\left( \ldots, \left(d_a\phi \right)_i, \ldots, \left(d_b\psi\right)_j, \ldots \right) 
 =-\left(\ldots, (d_b\phi)_i, \ldots,   
\left(d_a\psi\right)_j, \ldots\right).    
\end{eqnarray}
From $\phi$, $\psi \in \mathcal I(2)$ with $a=b$ it follows 
$d_a\phi$, $d_a\psi \in \mathcal I(2)$.  
Both with $d_b d_a =0$ or $d_b d_a \ne 0$, 
by applying further choices $d_c$, $c=a$ or $\overline{a}$ of differentials,  
we receive from \eqref{perdozo} further relations 
for higher differentials if the maximal powers 
$p\left(d_{a_n} \ldots d_{a_1}) \phi\right)$ 
of sequences 
$d_{a_n} \ldots d_{a_1}$, (where $a_i$, $1 \le i \le n$ 
is a choice of $a$ and $\overline{a}$) 
of the differentials $d_a$ and $d_{\overline{a}}$ permit. 
\subsection{Single differential, two kinds of elements}
Consider $d_a(\ldots, \phi^r, \ldots,  \psi^s, \ldots)$,  
with separate maximal powers for $\phi$ and for $\psi$,
 or with total maximal power $r+s=k$ for $\phi$ and $\psi$ together.  
It is assumed that $\phi$- and $\psi$-entries are mixed 
and commutation relations are not known. 
Let $d_a^{\alpha(\phi)} \phi=0$, $d_a^{\beta(\psi)}\psi=0$, 
where $\alpha(\phi)>1$ and $\beta(\psi)>1$  
are maximal orders of the differential $d_a$. 
Then, one has
\begin{eqnarray*}
  0=d_a \left(\ldots, \phi^r, \ldots,  \psi^s, \ldots \right) 
&=& \sum\limits_{i_1=1}^r  
\left(\ldots, (d_a\phi)_{i_1}, \ldots, (\phi)^{r-1}_{\widehat{i_1}},    
\ldots, (\psi)^s \ldots \right)
\nn
&+& \sum\limits_{i_2=1}^s  
\left(\ldots, (\phi)^r \ldots,   
(d_a\psi)_{i_2}, \ldots, (\psi)^{s-1}_{\widehat{i_2}}, \ldots \right).  
\end{eqnarray*}
By continuing the process of $n$-times application 
of the differential $d_a$ until  
 the powers of $\phi$ and $\psi$ become zero, 
  the sequence of identities stops. 
\subsection{Two differentials, identical elements of the total maximal power}   
For $q(\phi)$ identical elements $\phi$,
placed at all $i$-th positions,  
 $1 \le i \le q(\phi) \in \N$, such that 
 $\phi \in \mathcal I(q(\phi))$, i.e., 
$\phi^{q(\phi)}=0$, equivalently,   
$(\ldots$, $(\phi)_1$, $\ldots$, $(\phi)_i$, $\ldots$, $(\phi)_{q(\phi)}, \ldots)=0$,   
(i.e., $\phi$ is placed $q(\phi)$ times in various places in the product). 
 Then, 
\begin{eqnarray}
\label{normak}
&& 0=d_a\left( \ldots, (\phi)_i, 
\ldots \right)   
= \sum\limits_{i=1}^{q(\phi)} \left( \ldots, (d_a\phi)_i,  
\ldots, (\phi)_{\widehat{i}}^{q(\phi)-1},  \ldots \right), 
\end{eqnarray}
which is not trivial for $p(\phi) > 1$. 
In the case when $d_b d_a=0$, $b=a$ or $\overline{a}$, 
  by applying $d_b$ to \eqref{normak}, one finds   
\begin{eqnarray}
\label{normak1}
&& 0= \sum\limits_{i =1}^{q(\phi)} \sum\limits_{1 \le j \ne i}^{q(\phi)-1} 
\left(\ldots, \left(d_a\phi\right)_i, \ldots, \left(d_b \phi\right)_j, 
\ldots, (\phi)_{\widehat{i}, \widehat{j}}^{q(\phi)-2}, \ldots \right).       
\end{eqnarray}
With $(d_a\phi, d_b\phi) \in \mathcal I(2)$, 
this identity trivializes. 
For $b=a$ and  
$d_a\phi \notin \mathcal I(2)$,  
 one gets relations on 
differentials of $\phi$. 
Both for $d_b d_a = 0$ or $d_b d_a \ne 0$,  
by applying a sequence $d_{a_n} \ldots d_{a_1}$ 
of the differentials $d_a$ and $d_{\overline{a}}$ to \eqref{normak},  
the hierarchy of identities arises, 
\begin{eqnarray*}
0= \sum\limits_{{i_n, i_{n-1}, \ldots, i_1=1} 
\atop {\; i_s \ne j_r, \; 1 \le s \ne r \le n}}
^{q(\phi)-n, \ldots, q(\phi)-1, q(\phi)}  
\left( \ldots,  \left(d_{a_n} \ldots d_{a_1}\phi\right)_{i_1},  
\ldots, \left(d_{a_{n-1}} \ldots d_{a_1} \phi\right)_{i_n}, \ldots, 
(\phi)_{\widehat{i_1}, \ldots, \widehat{i_n}}^{q(\phi)-n}, \ldots 
 \right), &&
\end{eqnarray*}
when maximal orders and powers
 of corresponding combinations of differentials permit. 
\subsection{Orders of two differentials of several identical elements}
In this subsection we show how to use transfer of differentials 
for higher order differentials. 
Let $d_a^{p(\phi)}\phi=0$, and $\left(d_a^{p'}\phi\right)^{q(p', \phi)}=0$ 
for $p' < p(\phi)$ (we do not allow the ambiguity $0^0$). 
Then, for $0 < p_i < p(\phi)$,
  $1 \le q_i \le q(p_i, \phi)$, $1 \le i \le k$, 
 the identity is   
\begin{eqnarray} 
\label{porox}
&& 0=d_a \sum\limits_{i=1}^k \left(  \ldots, \left(d_a^{p_1} \phi\right)^{q_1},   
\ldots, \left(d_a^{p_i} \phi\right)^{q(p_i, \phi)}, \ldots, 
\left(d_a^{p_k} \phi\right)^{q_k}, \ldots \right) 
\nn
&&= \sum\limits_{i=1}^k \sum\limits_{s_1=1}^{q(p_i, \phi)}
 \left( \ldots, \left(d_a^{p_i+1} \phi\right)_{s_1},  
\ldots, \left(d_a^{p_i} \phi\right)^{q(p_i)-1}_{\widehat{s_1}},     
 \ldots, \left(d_a^{p_{\widehat{i}}} \phi\right)^{q_{\widehat{i}}}, 
 \ldots \right),  
\end{eqnarray}
which is non-trivial
if $p_i+1 < p(\phi)$, and $q(p_i+1, \phi) > 1$. 
Applying the differential $d_a$ further times  
one obtains higher identities. 
Then we can have a finite number $n$ applications of the differential 
 and corresponding sequence of identities.  
Similar products including also powers 
of the differential $d_{\overline{a}}$ 
lead to similar identities. 
\subsection{Two differentials, multiple elements}
Recall again that the main idea of differential identity generation 
is to be able to use them to transfer powers of $\mathcal C$-elements
into orders of corresponding differentials. 
Let us look at $(\ldots, (\phi_i)^{r_i}, \ldots)$, 
with distributed power of single elements  
 $\phi_i$-entries for all $i$, $1 \le i \le k \le l$, 
$1 < r_i \le q(\phi_i)$,  $(\phi_i)^{q(\phi_i)}=0$, i.e.,  
where $k$ is the total number of various $\phi_i$  
with individual maximal powers 
depending on an element $\phi_i \in \mathcal C$.    
Let $d_a^{p(\phi_i)} \phi_i=0$, $p(\phi_i) >1$, 
 $d_{\overline{a}}^{\overline{p}(\phi_i)} \phi_i=0$, 
$\overline{p}(\phi_i)>1$,  
and   
$\left(d_a^{p(\phi_i)} \phi_i\right)^{q(p_i, \phi_i)}=0$, 
$q(p_i, \phi_i) >1$,  
$\left(d_{\overline{a}}^{\overline{p}(\phi_i)} 
\phi_i\right)^{q(\overline{p}_i, \phi_i)}=0$, 
$q(\overline{p}_i, \phi_i) >1$, 
be the maximal orders and powers of the differentials
 $d_a$, $d_{\overline{a}}$ depending on $\phi_i$.  

 Consider the case of a single differential and multiple elements.  
In order to form the vanishing product, we 
include exactly one maximal power $q(\phi_i)$ 
of an element $\phi_i$, and powers 
$1 \le r_{\widehat{i}} < q(\phi_{\widehat{i}})$
 of other types of elements $\phi_{\widehat{i}}$, $1 \le {\widehat{i}} \ne i\le k-1$. 
For $d_a$ as before 
\begin{eqnarray}
 0&=& d_a\sum\limits_{i=1}^k \left( \ldots, 
\phi_1^{r_1},
 \ldots, \phi_i^{q(\phi_i)}, \ldots, (\phi_k)^{r_k}, \ldots \right)  
\nn
&=& \sum\limits_{i=1}^k  \sum\limits_{s=1}^{q(\phi_i)}  
\left(\ldots, (d_a \phi_i)_s, \ldots, (\phi_i)^{q(\phi_i)-1}_{\widehat{s}},  
 \ldots, \phi_{\widehat{i}}^{q(\phi_{\widehat{i}})},  \ldots \right).  
\end{eqnarray}
Further several applications of the differential $d_a$
 to the last identity results in further identities. 
It is clear that the orders of differentials 
$d_a(\phi_i)$ and $d_a(\phi_j)$ are growing until they reach $p(\phi_i)$ and 
$p(\phi_j)$, 
and thus corresponding summands vanish.  
Similarly, the powers of the differentials 
$d_a^{p_i} \phi_i$ and $d_a^{p_j} \phi_j$  
are growing until they reach $q(p_i, \phi_i)$ and $q(p_j, \phi_j)$, 
and corresponding summands vanish. 
\subsection{The most general polynomial case:
 the identities for orders of two differentials 
 and powers of multiple elements}   
In this subsection the identities for 
the most general 
polynomial differential-order and power products are computed. 
 Let us first assume that the maximal orders 
$p(\phi_i) > 1$ and $\overline{p}(\phi_i) > 1$  
of the differentials $d_a$ and $d_{\overline{a}}$ correspondingly, 
 i.e.,  such that $d_a^{p(\phi_i)}=0$, 
$d_{\overline{a}}^{\overline{p}(\phi_i)}=0$   
and the maximal powers $q(p_i, \phi_i) > 1$  
and $\overline{q}(\overline{p}_i, \phi_i) > 1$ of 
$d_a^{p_i(\phi_i)}$ and $d_{\overline{a}}^{\overline{p}_i(\phi_i)}$
i.e., such that 
$\left(d_a^{p_i(\phi_i)}\right)^{q(p_i, \phi_i)}=0$ 
and 
$\left(d_{\overline{a}}^{\overline{p}_i(\phi_i)}\right)
^{\overline{q}(\overline{p}_i, \phi_i)}=0$, 
do depend on elements $\phi_i \in \mathcal C$, 
for $1 \le i \le k$ types of elements $\phi_i$. 
 Take a product with exactly one differential in the maximal 
power $q(p_i, \phi_i)$, 
while all orders $p_i(\phi_i)$ of the differentials 
$d_a^{p_i(\phi_i)}$ are lower than the maximal order $p(\phi_i)$, 
i.e., $0 \le p_i(\phi_i) < p(\phi_i)$, 
\begin{eqnarray} 
\label{corex}
&& 0=d_a \sum\limits_{i=1}^k
 \left(\ldots, \left(d_a^{p_1(\phi_1)} \phi_1\right)^{r_1},  
\ldots,  \left(d_a^{p_i(\phi_i)} \phi_i\right)^{q_i(p_i, \phi_i)},     
\ldots, \left(d_a^{p_k(\phi_k)} \phi_k\right)^{r_k}, \ldots \right) 
\nn
&&= \sum\limits_{i=1}^k 
\sum\limits_{s_1=1}^{q(p_i, \phi_i)} \left( \ldots, (d_a^{p_i(\phi_i)+1} \phi_i)_{s_1},    
\ldots, (d_a^{p_i} \phi_i)^{q_i(p_i, \phi_i)-1}_{\widehat{s_1}},     
 \ldots, (d_a^{p_{\widehat{i}}} \phi_{\widehat{i}})^{r_{\widehat{i}}},  
 \ldots \right). 
\end{eqnarray}
It is clear that if we put at least one maximal 
order $p_i(\phi_i)=p(\phi_i)$ or $p_j(\phi_j)= p(\phi_j)$ 
  of the differential $d_a$ into the first line of \eqref{corex}
then all further identities trivialize. 

One can also take the case of the total maximal power elements. 
I.e., $\sum_{i=1}^t$ $q_i$ $(p_i, \phi_i)=q$, 
and $q_i(p_i, \phi_i) < q(p_i, \phi_i)$
 i.e., when the sum of powers reaches 
the critical value, and 
$\times\left(d_a^{p_i(\phi_i)} \phi \right)^{q_i(\phi_i)} \in \mathcal I(q)$.    
Similarly, $\sum_{i=1}^t p_i(\phi_i)=p$, and 
$p_i(\phi_i) < p(\phi_i)$,  
i.e., the sum of orders of the differentials 
is at its maximum, and 
$\times d_a^{p_i(\phi_i)} \phi_i \in \mathcal I(p)$
Then in both cases 
$\left( \ldots,  \left(d_a^{p_i(\phi_i)} \phi_i\right)^{q_i(p_i, \phi_i)},     
\ldots \right)=0$.
When we apply the second differential $d_{\overline{a}}$ to the identity 
\eqref{corex}, all branches of the 
hierarchy of identities arise 
 with appropriate conditions on orders and powers of differentials.  
\section{Examples of a single element closed products} 
In this Section we provide examples of closed products arising, in particular, 
from the identities of Section \ref{porgo}. 
Closed products and their identities are very useful in the theory of 
completely integrable \cite{arnold} and exactly solvable \cite{ls, razsav} systems. 
Note also possible applications for operads.  
Note that 
in the search of products annihilated by a differential, 
we take only non-vanishing products. 
To use the annihilation given by maximal powers of $d_a\phi_i$ 
in multi-element products after application of a differential $d_a$ once, 
we have to include not just one power of each element but 
lower powers also. 
The main idea of a closed product construction is  
in taking into account changes both in orders of differentials themselves 
and powers of their action on $\mathcal C$-elements. 
For that reason we have to include not only powers of elements 
acted by the differentials, but also powers of elements not acted by them. 
The balance of raising and lowering powers of elements results 
in extra conditions producing closed products. 
In some of the examples below we use the transfer of differentials 
using certain identities from Section \ref{porgo}. 
By integrating the identities of the form \eqref{porox} 
one finds sets of advanced closed products. 

Note again that when we write a product in the form 
$(\ldots, (\phi)_i, \ldots, (\psi)_j, \ldots)$, $1 \le i \ne j \le l$, 
it is not assumed that all elements $\psi$ are positioned on the right to 
all $\phi$ elements, i.e., $\phi$ and $\psi$ elements can be mixed in the product. 
\subsection{A single differential, a single element} 
First, let us look at expressions without differentials, 
$k < q(\phi)$ distributed position entries in a product, 
 where $q(\phi) >1$ is the maximal order of $\phi$, 
$\phi^{q(\phi)}=0$,
 i.e., $\phi \in \mathcal I(q(\phi))$.   
Let the maximal order of the differential $d_a$ is $p(\phi)$, 
i.e.,  $d^{p(\phi)}_a\phi=0$, 
and $q(p', \phi) \ge 1$ is the maximal power of $d^{p'}_a\phi$.  

 1.) With $p(\phi)=1$,  
\begin{eqnarray}
d_a\left(\ldots, (\phi)_1, \ldots, (\phi)_i, \ldots, (\phi)_k, \ldots\right)= 
\sum\limits_{i=1}^k 
\left(\ldots, (d_a\phi)_i, \ldots, (\phi)^{k-1}, \ldots \right)=0, 
\end{eqnarray}
 due to each summand is zero. 
 Since there is only one differential in each summand  
the use of the transfer formula \eqref{perdozo}
 would only changes the sign. 

2.) Suppose $d_ad_a \phi=0$, i.e., $p(\phi)=2$, and a distributed product   
 $(\ldots, (d_a\phi)_{q \; times, \ldots})=0$. 
We denote this relation as $(d_a\phi)^{q(1, \phi)}=0$, 
 i.e., $d_a\phi \in \mathcal I(q(1, \phi))$. 
 We continue with a distributed product containing $(d_a\phi)^k$   
with $k\le q(1, \phi)$ entries of $d_a\phi$. 
Then 
\begin{eqnarray}
d_a\left(\ldots, (d_a\phi)^k, \ldots\right)= \sum\limits_{i=1}^k 
\left(\ldots, (d_a d_a\phi)_i, \ldots, (d_a\phi)^{k-1}_{\widehat{i}}, \ldots\right)=0. 
\end{eqnarray}
3.) Next, with $d_ad_a\phi=0$, include $(\phi)^r$,   
i.e., take $(\ldots, (d_a\phi)^k, \ldots, (\phi)^r, \ldots)$, 
for $1 \le k+r \le l$, 
$k \le q(1, \phi)$ entries of $d_a\phi$, 
  $\left(d_a\phi\right)^{q(1, \phi)}=0$, for $d_a\phi$, 
 and $r \le q(\phi)$ entries of $\phi$, $\phi^{q(\phi)}=0$. 
Then, for $k+1=q(1, \phi)$,  
\begin{eqnarray*}
d_a\left(\ldots, (d_a\phi)^k, \ldots, (\phi)^r, \ldots\right)=  
\sum\limits_{s=1}^r \left(\ldots, (d_a\phi)^k, \ldots,
 (d_a\phi)_s, \ldots, (\phi)^{r-1}_{\widehat{s}}, \ldots\right)=0,  
\end{eqnarray*} 
 since the total power of $d_a\phi$ in each 
of summand becomes $q(1, \phi)$ and it, therefore, vanishes.  
\subsection{Higher orders of a single differential, single element} 
Now let us proceed with higher order closed products.
Let $p(\phi)$ be the maximal order of the differential $d_a\phi$, 
i.e., $d_a^{p(\phi)}\phi=0$.    

 1.) With $k$ entries of $d_a^{p_i}\phi$ differentials, 
$k < q(p_i, \phi)$, $\left(d_a^{p_i}\phi\right)^{q(p_i, \phi)}=0$, 
$1 \le p_i < p(\phi)$,  i.e., one gets the product  
$\left(\ldots, \left(d^{p_i}\phi\right)_i, \ldots \right)$. 
Then, for $p_i = p(\phi)-1$ for all $i$, 
\begin{eqnarray}
d_a\left(\ldots, \left(d_a^{p_i}\phi\right)_i, \ldots \right)= 
\sum\limits_{s=1}^k \left(\ldots, 
\left(d_a^{p(\phi)}\phi\right)_s, \ldots, 
\left(d_a^{p_i}\phi\right)_{\widehat{s}}, \ldots \right)=0.  
\end{eqnarray}

2.) Now, let $\phi^{q(\phi)}=0$.  Include 
$k < q(p_i, \phi)$,  entries of $d_a^{p_i}\phi$, $1 \le i \le k$, 
 and $1 \le r < q(\phi)$ entries of $\phi$, i.e., 
 $\left(\ldots, \left(d_a^{p_i}\phi\right)_i, \ldots, (\phi)^r, \ldots\right)$.  
One has
\begin{eqnarray}
\label{pordon}
&& d_a\left(\ldots, (d_a^{p_i}\phi)_i, \ldots, (\phi)^r, \ldots\right)= 
\sum\limits_{s=1}^k \left(\ldots, (d_a^{p_s+1}\phi)_s, \ldots, 
 \left(d_a^{p_{\widehat{s}}}\phi\right)_{\widehat{s}}, 
\ldots, (\phi)^r, \ldots \right)  
\nn
&&\qquad \qquad + \sum\limits_{t=1}^r 
 \left(\ldots, (d_a^{p_i}\phi)_i, \ldots, (d_a\phi)_t, \ldots, (\phi)^{r-1}, \ldots\right). 
\end{eqnarray}
For $p_i=1$, $k+1=q(1, \phi)$, and 
  $p_i+1=2=p(\phi)$, $1 \le s \le k$, both 
 groups of summands vanish. 
 Note that that does not depend on $r$. 
By inclusion of an extra summation into \eqref{pordon} 
one can use \eqref{normak1} to move differentials among $d_a\phi$ and $\phi$ 
with changing the sign. 
The symmetry of allows to distribute the orders
 of differentials among entries of $\phi$. 
\subsection{Orders of a single differential, single element, powers of differentials}
 1.) Let  
$\left(\ldots, \left(d_a^{p_i}\phi\right)^{q_i}_i, \ldots\right)$ be a product 
containing $k$ various powers of the differential $d_a$, $1 \le i \le k \le l$, 
for $0 \le p_i \le p(\phi)$, $1 \le q_i \le q(p_i, \phi)$,
 where $d_a^{p(\phi)}\phi=0$, 
$(d_a^{p_i} \phi)^{q(p_i, \phi)}=0$. 
We then obtain 
\begin{eqnarray*}
d_a\left(\ldots, \left(d_a^{p_i}\phi\right)^{q_i}_i, \ldots \right) &=& 
\sum\limits_{s=1}^k  \sum\limits_{t=1}^{q_s}  
\left(\ldots, \left(d_a^{p_s+1}\phi\right)_t, \ldots,  
\left(d_a^{p_s}\phi\right)^{q_s-1}_{\widehat{t}}, \ldots, 
  \left(d_a^{p_{\widehat{s}}}\phi\right)^{q_{\widehat{s}}}_{\widehat{s}},  \ldots\right), 
\end{eqnarray*}
where $\widehat{t}$ denote the omission of $t$. 
That product is closed when   
 1.) $p_i+1=p(\phi)$ for all $1 \le i \le k$; 
2.) $d_a^{p_s+1}\phi$, $d_a^{p_s}\phi \in \mathcal I(2)$; 
 3.) $d_a^{p_s+1}\phi$, $d_a^{p_s}\phi \in \mathcal I(q_s)$; 
4.) $d_a^{p_s+1}\phi$, $(d_a^{p_{\widehat{s}}} \phi) \in \mathcal I(2)$; 
5.) $d_a^{p_s+1}\phi$,
$(d_a^{p_{\widehat{s} } } \phi)^{q_{\widehat{s} } }$    
  $\in \mathcal I (q_{\widehat{s}}+1)$.  

2.) The general expression for a closed product 
in the case of a single element with 
 $k$ types of orders and powers of 
 a single differential is given by  
$\left(\ldots, \left(d_a^{p_i}\phi\right)^{q_i}_i, \ldots, 
\left(d_a^{p_i-1}\phi\right)^{r_i}_j \right.$,  $ \left.\ldots \right)$, 
where for any pair $p_i=p_j-1$, $2p_i< p_i(\phi)$, and 
$q_i+r_j< q(p_i, \phi)$, $1 \le i, j \le k$ 
(such that the initial product does not vanish). 
\begin{eqnarray*}
&& d_a\left(\ldots, \left(d_a^{p_i}\phi\right)^{q_i}_i, \ldots, 
\left(d_a^{p_i-1}\phi \right)^{r_i}_j, \ldots \right)   
\nn
&&
=\sum\limits_{s=1}^{k} \sum\limits_{t=1}^{q_s}  
 (\ldots, (d_a^{p_s+1}\phi)_t, \ldots, (d_a^{p_s } \phi)^{q_s-1}_{\widehat{t}}, 
\ldots, (d_a^{p_i-1}\phi)^{r_j}_j, \ldots ) 
\nn
&&
  + \sum\limits_{s=1}^k \sum\limits_{t=1}^{r_s}  (\ldots, (d_a^{p_i}\phi)^{q_i}_i, 
\ldots, (d_a^{p_s}\phi)_t, \ldots, 
(d_a^{p_s-1}\phi)_{\widehat{t}}^{r_s-1}, \ldots )=0, 
\end{eqnarray*}
when 1.) $p_i+1=p(\phi)$ and $q_i+1=q(p_i, \phi)$, 
 for all $1 \le i \le k$; 2.) $p_i+1=p(\phi)$ and $p_i=p_s$, $2p_i \ge p_i(\phi)$. 
\subsection{Two differentials, single element}
For the product $\left(\ldots, d_{\overline{a}}\phi, 
\ldots, d_a\phi, \ldots, (\phi)^r, \ldots\right)$, with $1 < r < q(\phi)$ one has: 
1.) For $d_a\phi \in \mathcal I(2)$, i.e., $(d_a\phi)^2=0$, and 
 $d_a d_{\overline{a}}\phi=0$, $d^2_a\phi=0$ 
\begin{eqnarray*}
 d_a\left(\ldots, d_{\overline{a}}\phi, \ldots, d_a\phi, \ldots, (\phi)^r, \ldots \right)  
=(\ldots, d_a (d_{\overline{a}}\phi), \ldots, d_a\phi, \ldots, (\phi)^r, \ldots)
 \qquad \qquad \qquad &&
\nn
 + (\ldots, d_{\overline{a}}\phi, \ldots, d_a d_a\phi, \ldots, (\phi)^r, \ldots)
 + (\ldots, d_{\overline{a}}\phi, \ldots, d_a\phi,
 \ldots, d_a\phi, \ldots, (\phi)^{r-1}, \ldots)=0. &&
\end{eqnarray*} 

2.) For $d_{\overline{a}} d_a \phi\ne 0$, $d_a d_{\overline{a}} d_a\phi= 0$, 
  $d_a^2\phi=0$, and $d_a\phi \in \mathcal I(2)$,  
the product vanishes: 
\begin{eqnarray*} 
  d_a\left( \ldots, d_{\overline{a}} d_a \phi,     
\ldots, d_a \phi, \ldots, (\phi)^r, \ldots \right) 
 = 
\left(\ldots, d_a d_{\overline{a}} d_a \phi,
 \ldots, d_a\phi, \ldots,  (\phi)^r, \ldots\right) \qquad \qquad \qquad  &&
\nn
+ \left( \ldots, d_{\overline{a}} d_a \phi,     
\ldots, d_ad_a \phi, \ldots, (\phi)^r, \ldots \right) 
+ \left( \ldots, d_{\overline{a}} d_a \phi,     
\ldots, d_a \phi, \ldots, d_a\phi, \ldots, (\phi)^{r-1}, \ldots\right).\qquad &&
\end{eqnarray*}
Note that both closed products do not depend on $r$.

When commutation rules for elements in a product are known,  
we are able to move them around. 
In a distributed case, we 
first gather the powers of the same element together.  
Then the classical formula 
$d_b\left(d_a^p \phi\right)^q= q d_b\left(d_a^p \phi\right) 
\left(d_a^p \phi\right)^{q-1}$ 
(we agree to put the extra derivative in front of the power of the derivative) 
for a differential applied to a power 
of a $\mathcal C$-element is true, and we can write invariant elements 
in the explicit form. 
\section{Examples of multi-element closed products}
\label{mima} 
In this Section we give examples of multi-element closed product.  
\subsection{Two-elements closed product}
Here we have the case of dependence on two elements $\phi$ and $\psi$. 
That case does not fall in the general idea \cite{Ghys, Ko}
of one-parameter-element 
invariants. 
Let us see if it is possible to transpose $d_a$-differential. 
1.) For $(d_a\phi$, $d_a\psi) \in \mathcal I(2)$, and $d_ad_a \phi=0$,  
\begin{eqnarray*}
 d_a(\ldots, d_a\phi, \ldots, \psi, \ldots)=
  (\ldots, d_ad_a\phi, \ldots, \psi, \ldots)  
 + (\ldots, d_a\phi, \ldots, d_a\psi, \ldots)=0. 
\end{eqnarray*}

2.) With $(d_a\phi$, $d_{\overline{a}}\psi) \in \mathcal I(2)$, 
 $d_a d_{\overline{a}}=0$, then   
\begin{eqnarray*}
 d_a(\ldots, d_{\overline{a}}\phi, \ldots, \psi, \ldots)  
 =(\ldots, d_a d_{\overline{a}} \phi, \ldots, \psi, \ldots) 
 + (\ldots, d_a\phi, \ldots, d_{\overline{a}} \psi, \ldots)=0;  
\end{eqnarray*}

 3.) $d_a\phi$, $\psi \in \mathcal I(2)$,
 $d_a d_{\overline{a}}\phi=d_{\overline{a}} d_a\phi=0$,  
\begin{eqnarray*}
 d_a(\ldots, d_{\overline{a}} \phi, \ldots, \psi, \ldots)  
 &=&(\ldots, d_{\overline{a}} d_a\phi, \ldots, \psi, \ldots) 
 + (\ldots, d_a\phi, \ldots, d_{\overline{a}}\psi, \ldots) 
\nn
 &=&-(\ldots, d_a\phi, \ldots, d_{\overline{a}}\psi, \ldots)  
 + (\ldots, d_a\phi, \ldots, d_{\overline{a}}\psi, \ldots)=0. 
\end{eqnarray*}
due to property \eqref{perdozo}. 
\subsection{Two differentials, three elements}
Take the product 
 $(\ldots, d_{\overline{a}}\phi, \ldots, d_a\psi, \ldots, (\chi)^r, \ldots)$ with 
 $d_ad_a\psi=0$, $d_ad_{\overline{a}}\phi=0$, and   
1.) $d_{\overline{a}} \phi$, $d_a\psi \in \mathcal I(2)$; 
2.) $d_{\overline{a}} \phi$, $d_a\chi \in \mathcal I(2)$;
3.)   $d_a\psi$, $d_a\chi \in \mathcal I(2)$;   
4.) $\psi$, $d_a\chi \in \mathcal I(2)$, $d_ad_a\chi=0$; 
 5.) $\phi$, $d_a\chi \in \mathcal I(2)$, and $d_{\overline{a}} d_a\chi=0$, 
6.) $\phi$, $d_a\psi \in \mathcal I(2)$, and $d_{\overline{a}} d_a\psi=0$. 
Then the product vanishes
\begin{eqnarray*}
 d_a(\ldots, d_{\overline{a}}\phi, \ldots, d_a\psi, \ldots, (\chi)^r, \ldots)   
= (\ldots, d_a(d_{\overline{a}}\phi), \ldots, d_a\psi, \ldots, (\chi)^r, \ldots) &&
\nn
+ (\ldots, d_{\overline{a}}\phi, \ldots, d_a d_a\psi, \ldots, (\chi)^r, \ldots)  
+ (\ldots, d_{\overline{a}}\phi, \ldots, d_a\psi, \ldots, d_a\chi, 
\ldots, (\chi^{r-1}), \ldots), &&
\end{eqnarray*}
according to property \eqref{perdozo}.
\subsection{Orders of a single differential, powers of several elements} 
Consider the following case: 
$\left(\ldots, (d_a^{p_i}\phi_i)^{q_i(p_i)}, \ldots,   
(\phi_j)^{q_j},  \ldots \right)$, 
with $0 \le  p_i(\phi_i) < p(\phi_i)$, and 
$1 \le  q_i(p_i) < q(p_i, \phi_i)$, $1 \le  q_j(p_j) < q(p_j, \phi_j)$,
$1 \le  q_i < q(\phi_i)$, $1 \le  q_j < q(\phi_j)$, 
$1 \le i \le k$,  $1 \le j \le k'$,  
 the maximal order and power 
(here it starts from $1$ to preserve $l$, 
or keep the identical element in the product).   
Note that the powers are lower than 
  $p(\phi_i)-1$ and $q(\phi_j)-1$ correspondingly 
 to insure that the corresponding product does not vanish. 
Next, by applying conditions on $p_i(\phi_i)$ and $q(\phi_i)$, 
we find their combinations such that the product is closed.  

1.) It is clear that for all $p_i+1=p(\phi_i)$, $1 \le i \le k$,  
$p_i=1$, and $ q_i(1)+1=q(1)$;  
2.) for all $p_i+1=p(\phi_i)$, $1 \le i \le k$, and 
with $p_i\ne 1$ , $d_a\phi_i$, $\phi_j \in \mathcal I(2)$;   
3.) for all $p_i+1=p(\phi_i)$, $1 \le i \le k$, and 
with $p_i\ne 1$ , $d_a\phi_i$, $(\phi_j)^{q_j-1} \in \mathcal I(q_j)$ 
the product is closed. 
\begin{eqnarray*}
&& d_a \left(\ldots, (d_a^{p_i}\phi_i)^{q_i(p_i)}, \ldots,   
(\phi_i)^{q_i},  \ldots \right)
\nn
&&
=\sum\limits_{s=1}^k \sum\limits_{t=1}^{q_s(p_s)} 
\left(\ldots, (d_a^{p_s+1}\phi_s)_t, \ldots, 
 (d_a^{p_s}\phi_s)^{q(p_s)-1}_{\widehat{t}}, \ldots,     
(\phi_i)^{q_i},  \ldots \right) 
\nn
&&
+ \sum\limits_{s'=1}^{k'} \sum\limits_{t=1}^{q_{s'}} 
\left(\ldots, 
 (d_a^{p_i}\phi_i)^{q_i(p_i)}, \ldots,     
(d_a \phi_{s'})_t, \ldots,  (\phi_{s'})^{q_{s'}-1}_{\widehat{t}},  \ldots \right), 
\end{eqnarray*}
where $\widehat{t}$ denotes the omission of the index $t$. 
\subsection{Powers of orders of multiple differentials, multiple elements}
Consider the product
$\left(\ldots, (d^{p_{n, i}}_{a_{n, i}} \ldots d^{p_{1, i}}_{a_{1, i}}\phi_i)_i
^{q({\bf p}_{n, i}, {\bf a}_{n, i})}, \ldots \right)$, 
for some $1 \le i \le k \le l$, and where we use the notation 
 ${\bf x}_{n, i}= (x_{n, i}, x_{n-1, i}, \ldots, x_{1, i})$. 
For the product not to vanish, all orders $p_{s, i}$, $1 \le s \le n$ satisfy 
$p_{s, i}\left(d^{p_{s-1, i}}_{a_{s-1, i}} \ldots d^{s, i}_{a_{1, i}} \phi\right) 
< p\left(d^{p_{s-1, i}}_{a_{s-1, i}} \ldots d^{p_1, i}_{a_{1, i}} \phi \right)$, 
where 
$p\left(d^{p_{s-1, i}}_{a_{s-1, i}} \ldots d^{p_1, i}_{a_{1, i}} \phi\right)$ is the 
maximal order of the corresponding multiple differentials, 
and the power of the corresponding multiple differentials 
$d^{p_{n, i}}_{a_{n, i}} \ldots d^{p_{1, i}}_{a_{1, i}}\phi_i$ is lower than 
${q({\bf p}_{n, i}, {\bf a}_{n, i})} < q({\bf p}_{n, i}, {\bf a}_{n, i}, \phi)$
is lower than the maximal ones for all $i$. 
Then 
\begin{eqnarray}
d^{p_{n+1}}_{a_{n+1}}\left(\ldots, 
\left(d^{p_{n, i}}_{a_{n, i}} \ldots d^{p_{1, i}}_{a_{1, i}}\phi_i\right)_i  
^{q({\bf p}_{n, i}, {\bf a}_{n, i})}, \ldots\right)=0,  
\end{eqnarray}
when 
1.) there exists $i$, $1 \le i \le k$, such that 
 $p_{n+1, i}\left(d^{p_{n, i}}_{a_{n, i}} \right.$ $\ldots $  
 $d^{p_1, i}_{a_{1, i}}$    
 $ \left. \phi\right)$
 $\ge$  
 $p\left(d^{p_{n+1, i}}_{a_{n, i}} \ldots d^{p_1, i}_{a_{1, i}} \phi\right)$; 
2.) if $k \ge 1$, 
and  $({\bf p}_{n+1, i}, {\bf a}_{n+1, i})= ({\bf p}_{n, j}, {\bf a}_{n, j})$, 
for some $1 \le i \ne j \le k$. 
\section{Multiple-graded differential algebras}
\label{bardos}
Let $\gamma_s \in \mathcal C$, $s \in \Z$, and      
${\bm J}_s=\left(J_{s, 1}, \ldots,  J_{s, n_s}\right)$ be a set   
 of choices $J_{s, j}=\left({a_{s, j} \atop r_{s, j}}\right)$, 
$1 \le j \le n_s$. 
Let us require that for $l$ chain-cochain spaces  
$C\left(\Theta^{{\bm n_i}}_{\bm m_i}\right)$ 
of the multiple complex $\mathcal C$,    
 there exist subspaces 
$C\left(\Theta^{'\bm {n_i}}_{\bm m_i}\right) 
\subset C\left(\Theta^{{\bm n_i}}_{\bm m_i}\right)$    
 such that for all 
$\phi_{s, i} \in C\left(\Theta^{'{\bm n}_{s, i}}_{{\bm m}_{s, i}}\right)$,     
$1 \le i \le k$, 
\begin{equation}
\label{groby}
 \mathcal D_{{\bm J}_s} \gamma_s
=\cdot_l \left(\ldots, 
\left(\mathcal D_{\bm J_{s, i}} \phi_{s, i}\right)^{q_{s, i}}, \ldots \right),  
\end{equation} 
where it is clear that 
$q_{s,i} < q\left(\mathcal D_{\bm J_{s, i}} \phi_{s, i}\right)$
 are less than the maximal powers of the corresponding elements 
$\mathcal D_{\bm J_{s, i}} \phi_{s, i}$, 
and the maximal orders of differentials 
constituting (as in \eqref{pirdos}) $\mathcal D_{\bm J_{s, i}} \phi_{s, i}$  
satisfy the recurrence conditions \eqref{pirdo}.  
 We call this a differential condition. 
Symmetrizing with respect to all choices of ${\bm J_s}$ for a fixed $l$,  
we obtain a set of differential conditions 
\begin{equation}
\label{mordo}
\left\{  Symm_{{\bm J}_s, l} 
\left\{ \mathcal D_{\bm J_s} \gamma_s=   
 \; \cdot_l \left(\ldots, \left(\mathcal D_{\bm J_{s, i}}     
\phi_{s, i}\right)^{q(i, s)}, \ldots \right) \right\} \right\}, 
\end{equation}
for $\gamma_s \in \mathcal C$, $s \in \Z$. 
With $\mathcal D_{\bm J_s} \gamma_s=0$ for some $s$, 
 we call \eqref{mordo} an orthogonality condition. 
For each differential condition above of the set \eqref{mordo} 
  we have the coherence condition for 
corresponding upper and lower indices of $\mathcal C$-spaces, i.e., 
${\bm n}$ $\left( \mathcal D_{J_s} \gamma_s\right)$    
 $=$ $\sum_{i=1}^l \left( {\bm n_i} \left(\mathcal D_{\bm J_{s, i}} \phi_{s, i} \right)  
 - r_i\left(\mathcal D_{\bm J_{s, i}} \phi_{s, i}\right) \right)$,   
${\bm m}\left(\mathcal D_{\bm J_s} \gamma_s\right)  
=\sum_{i=1}^l 
\left({\bm m}_i\left(\mathcal D_{\bm J_{s, i}} \phi_{s, i}\right) \right.$ 
$ \left. - {\bm t}_i\left(\mathcal D_{\bm J_{s, i}} \phi_{s, i} \right) \right)$.    
We will skip these coherence relations  
for all differential conditions below. 
The notion of a set of differential and orthogonality
 conditions \eqref{mordo} generalizes  
corresponding orthogonality conditions assumed in \cite{Ghys, Ko}. 
One could like to understand which part 
of differential conditions \eqref{mordo} 
are independent. 
Indeed, in the case with known commutation rules for 
$\mathcal C$-elements and differentials,  
we use \eqref{roda} to normalize its sequence in the definition \eqref{pirdos} 
and in the conditions \eqref{mordo} 
by sending all $d_{\overline{a}}$ to the left with respect to $d_a$. 
Then we use known vanishing rules for powers of differentials. 
In the case when commutation rules are not known, 
one might use the independence of powers of differentials. 
The next result of this paper is the following. 
\begin{lemma}
\label{kuz} 
The rules \eqref{roda}, the differential conditions \eqref{mordo},   
and a set of maximal orders and powers for all $\mathcal C$-elements    
endow $\mathcal C$      
with the structure of a multiply graded  
infinite-dimensional
 differential algebra 
with respect to a $\cdot_l$-multiplication, $l \ge 0$. 
\end{lemma}
Note that, apart from the value of $l$ for a multiple product, 
the parameters characterizing a differential algebra above 
are the distributions of the net of mixed or non-mixed ideals 
$\mathcal I(p)$, $\mathcal I(q) \in \mathcal C$; 
the distributions of orders and powers of differentials 
 bounded by their maximal values when applied 
to $\mathcal C$-elements; domains of values of indices for $\mathcal C$-spaces; 
the set of rules \eqref{roda} for differentials; 
 conditions on completions with respect to differentials; 
commutation relations for $\mathcal C$-elements (when known); 
the domains of values for $\Theta^{\bm n}_{\bm m}$ parameters for 
$C\left(\Theta^{\bm n}_{\bm m}\right)$, and distribution of horizontal and vertical 
indices for the corresponding differentials.  
Now let us give a proof of the Lemma.
\begin{proof} 
The way of operation with sequences 
of differential and orthogonality conditions
 follows from 
the rules \eqref{roda}, Leibniz rule \eqref{leibniz}, 
 as well as from taking intro account maximal orders and powers of 
differentials and elements, 
i.e., 
the assumption that some of $\mathcal C$ belong to the ideals 
with respect to orders 
$\mathcal I(p)$ or powers $\mathcal I(q)$.  
Starting from a particular orthogonality or differential condition, 
we act consequently by the differentials $d_a$ and $d_{\overline{a}}$. 
Secondly, using the $\mathcal I(p)$, $\mathcal I(q)$ 
 ideal vanishing properties 
applied to orthogonality conditions, 
and by using the completeness property of  
the complex $\mathcal C$ with respect to the $\cdot_l$-product, 
we express particular elements $d^{j_i} \phi_i$ 
in terms of other $\mathcal C$ elements. 
With the maximal orders and powers for each 
particular $\mathcal C$-element,   
 a differential or an $\mathcal C$-element  
reaches its maximal order or power, 
and corresponding summand vanishes. 
Continuing the process,
 we finally obtain the full structure of differential conditions.  
A sequence of relations does not not stop 
as long as coherence conditions on indices are fulfilled, 
of until the sequence gives the both side identical zero. 
In some cases, due to the completeness condition for the $\cdot_l$-product,  
a sequence of differential conditions becomes infinite in a certain branch 
of the hierarchy. 
Let us reproduce the general structure of relations 
following from differential conditions for 
 an expression in the form of 
a multiple product \eqref{tugo}.
 The most general configuration associated 
to a $\cdot_l$-multiple product where several $\mathcal C$-elements is with elements that 
are situated at various places and mixed with the corresponding completions. 
In our setup, we work with  
 differential conditions containing powers of elements of $\mathcal C$
as well as orders of the differentials. 
The general form of the differential condition is given by \eqref{mordo},   
where $\mathcal D_{\bm J_{s, i}}\phi_{s, i}$
 is of the form \eqref{pirdos},  
$1 \le i \le k \le l$. 
By applying the differentials to \eqref{mordo} one arrives at further 
differential and orthogonality relations.  
Recall, that for each differential condition in this proof, 
there exist a coherence relation for corresponding $\mathcal C$-indices. 
The sequences of differential and orthogonality conditions 
derived above 
together with corresponding coherences relations for indices, 
taken for all  
choices of ${\bm J}_s$, provide the full set of differential 
and orthogonality relations defining 
the multiply graded differential algebra. 
The structure of a resulting differential algebra relations
as well as the structure of corresponding closed products 
given by Theorem \ref{buldo} 
 depends on which differential or orthogonality condition 
of a hierarchy we started from. 
In practice,  
one can start with a particular set 
of differential or/and orthogonality conditions. 
Then the resulting multiple graded differential algebra 
is a reduction of the full algebra. 
In addition to that, one can also restrict 
 domain of definitions for gradings for some families of 
chain-cochain complexes. 
E.g., one can set $n\ge 0$ instead of $\Z$. 
This will affect the coherence conditions for corresponding identities. 
\end{proof}
\section*{Acknowledgment}
The second author is supported by 
the Institute of Mathematics, Academy of Sciences of the Czech Republic (RVO 67985840).  

\end{document}